\documentclass[preprint,3p,11pt]{elsarticle}

 \usepackage{hyperref}
 \hypersetup{backref, colorlinks=true, linkcolor=blue}

\usepackage{stfloats}
\usepackage{subcaption}

\usepackage{graphicx} % allows inclusion of graphics
\usepackage{booktabs} % nice rules (thick lines) for tables
\usepackage{microtype} % improves typography for PDF
\usepackage{amssymb,amsmath}
\usepackage[ruled,vlined]{algorithm2e}
\usepackage{xcolor}
\usepackage[normalem]{ulem}

\graphicspath{{./figs/}}

\begin{document}
\begin{center}
{\bf \Large Analysis of Hybrid MC/Deterministic Methods for Transport Problems 
\\
\medskip
 Based on Low-Order Equations Discretized by Finite Volume Schemes}
\end{center}
\begin{center}
Vincent N. Novellino and Dmitriy Y. Anistratov 
\end{center}
\begin{center}
{ \it Department of Nuclear Engineering,
North Carolina State University,
Raleigh, NC 27695 \\
 vnnovell@ncsu.edu, anistratov@ncsu.edu }
\end{center}

\begin{frontmatter}
\begin{abstract}
 This paper presents
hybrid numerical techniques for solving the Boltzmann transport equation formulated by means of  low-order equations 
for angular moments of the angular flux.
The moment equations   are derived by the projection operator approach. 
The projected equations are closed exactly  using a high-order transport solution.
The low-order  equations of the hybrid methods are approximated with a finite volume   scheme of the second-order accuracy.
Functionals defining the closures in the discretized low-order equations are calculated  by Monte Carlo techniques.
In this study, we analyze effects of statistical noise and discretization error
on the accuracy of the hybrid transport solution.  
\end{abstract}
\end{frontmatter}

\section{Introduction}

Monte Carlo (MC) and deterministic methods for solving neutron transport problems are widely used by modeling software to make informed design decisions.
MC is generally slow to converge when seeking a global solution.
Hybrid MC/deterministic methods have been developed for solving particle transport problems to improve efficiency of MC simulations.
They can reduce statistical variance of the solution and accelerate MC calculations.
A group of hybrid methods for global transport problems are based on different low-order equations for moments of the angular flux \cite{larsen-yang-nse-2008,lee-physor-2010,wolters-nse-2013,pozulp-mc2023}. 

We develop hybrid numerical techniques applying low-order equations of the Quasidiffusion (QD) and Second Moment (SM) methods \cite{gol'din-cmmp-1964,sm-1976,mla-ewl-pne-2002}. 
Low-order equations for moments of the angular flux were derived by the projection operator approach. 
Exact closures were formulated for the projected equations using a high-order transport solution.
The low-order moment equations of the hybrid methods were discretized with a finite volume (FV) scheme.
Functionals defining the closures in the discrete moment equations are computed by MC.
There are two main factors affecting the numerical solution of hybrid methods: (1) the discretization error related to the finite size of mesh cells and (2) the statistical variation of functionals associated with finite number of particle histories.
In this study, we analyze effects of statistical noise and discretization error on the accuracy of the hybrid transport solution.

\section{Hybrid Transport Methods}

We consider one-group steady-state  transport problem in 1D slab geometry with isotropic scattering and source:
\begin{equation} \label{t-eq}
  \mu \frac{\partial \psi}{\partial x}(x,\mu) + \Sigma_t(x) \psi(x,\mu) = \frac{\Sigma_s(x)}{2} \int_{-1}^1 \psi(x,\mu') d\mu' + \frac{Q(x)}{2} \, 
\end{equation}
\begin{equation}
  \psi(0,\mu) = 0 , \ \mu > 0 \, , \quad 
  \psi(L,\mu) = 0 , \ \mu < 0 \, .
 \end{equation}
The low-order QD (LOQD) equation in the second-order form is given by \cite{gol'din-cmmp-1964}
\begin{equation}
-\frac{d}{dx}   \frac{1}{\Sigma_t(x)} \frac{d(E\phi)}{dx}(x)  + \Sigma_a(x) \phi(x) = Q(x) \, ,
\end{equation}
where the QD (Eddington) factor is defined by
\begin{equation} \label{QDf}
 E(x) = \frac{\int_{-1}^{1} \mu^2 \psi(x,\mu) d \mu}{  \int_{-1}^{1} \psi(x,\mu) d \mu} \, .
\end{equation}
The low-order SM (LOSM) equation has the  form  \cite{sm-1976}:
\begin{equation}
-\frac{d}{dx}   \frac{1}{3\Sigma_t(x)} \frac{d\phi}{dx}(x)  + \Sigma_a(x) \phi(x) =
 Q(x) +\frac{d}{dx}   \frac{1}{3\Sigma_t(x)} \frac{dF}{dx}(x)  \, ,
 \vspace{-0.15cm}
\end{equation}
where the SM factor is given by
\begin{equation}\label{SMf}
 F(x) = \int_{-1}^{1} \bigg(\frac{1}{3} -  \mu^2 \bigg)\psi(x,\mu) d \mu \, .
\end{equation}
The details on formulation of boundary conditions (BCs) for LOQD and LOSM equations are omitted for sake of brevity (see  \cite{gol'din-cmmp-1964,mla-ewl-pne-2002}).
QD and SM factors (Eqs. \eqref{QDf} and \eqref{SMf}) define the closures for the low-order equations and the high-order transport equation (Eq. \eqref{t-eq}). 

We formulate the hybrid QD (HQD) and hybrid  SM (HSM) methods for calculation of the scalar flux as a solution of the discrete low-order equations with the factors $E$ and $F$, respectively, computed with MC.
The low-order equations are discretized with a finite volume (FV) scheme of the second-order accuracy.
The discrete FV LOQD equations for the cell-average scalar flux $\phi_i$ have the following form:
\begin{equation}
-  \frac{ E_{i+1}\phi_{i+1} - E_i \phi_i  }{\Sigma_{t,i+1/2} \Delta x_{i+1/2}} 
+ \frac{ E_i\phi_i - E_{i-1} \phi_{i-1}  }{\Sigma_{t,i-1/2} \Delta x_{i-1/2}}    +
\Sigma_{a,i} \Delta x_i \phi_i =  Q_i \Delta x_i \, ,
\end{equation}
where $i$ is the spatial cell index, the integer subscript indicates a quantity averaged over the $i^{th}$ cell $[x_i, x_{i+1}]$,  $\Delta x_i = x_{i+1} - x_i$,
\begin{equation}
\Sigma_{t,i+1/2} =\frac{\Sigma_{t,i}\Delta x_i  + \Sigma_{t,i+1}\Delta x_{i+1}}{\Delta x_i  + \Delta x_{i+1}} \, , \quad
\Delta x_{i+1/2}=\frac{1}{2}(\Delta x_i+\Delta x_{i+1}) \, .
\end{equation}
The FV scheme for the LOSM equation is given by
\begin{equation}
-  \frac{ \phi_{i+1} - \phi_i  }{3\Sigma_{t,i+1/2} \Delta x_{i+1/2}} 
+ \frac{ \phi_i - \phi_{i-1}  }{3\Sigma_{t,i-1/2} \Delta x_{i-1/2}}    +
\Sigma_{a,i} \Delta x_i \phi_i =  Q_i \Delta x_i  
 + \frac{F_{i+1} - F_i  }{\Sigma_{t,i+1/2} \Delta x_{i+1/2}} 
+ \frac{ F_i  - F_{i-1} }{\Sigma_{t,i-1/2} \Delta x_{i-1/2}} \,.
\end{equation}

Track-length tallies are used for the functional estimators to compute cell-average factors
\begin{equation} 
  E_i    =   \frac{\sum_{n=1}^N   \sum_{m=1}^{M_n} \mu_{n,m}^2  w_{n,m} \ell_{n,m}}
              {\sum_{n=1}^N \! \! \sum_{m=1}^{M_n}   w_{n,m} \ell_{n,m}} \, ,
\end{equation}
\begin{equation} 
   F_i     =    \frac{1}{N\Delta x_i}\sum_{n=1}^N  \sum_{m=1}^{M_n}   \Big(\frac{1}{3} -  \mu_{n,m}^2 \Big) w_{n,m} \ell_{n,m} \, ,
\end{equation}
where $N$ is the number of particle histories, $M_n$ is the number of tracks of $n^{th}$ particle in the $i^{th}$ cell, $\ell_{n,m}$ is the $m^{th}$ track length of the $n^{th}$ particle in the $i^{th}$ cell, $w_{n,m}$ is the weight of this particle.
A general algorithm of the HQD and HSM  methods is described by Algorithm \ref*{alg:hybrid_algorithm}.
\begin{algorithm}
  \DontPrintSemicolon
  %\dontprintsemicolon
  \ForEach{$n=1\dots N$}{
     Simulate particle history and collect moments of tracklength and facet crossing tallies\;
  }
  Calculate functionals\;
  Solve the low-order equations  to compute $\phi_i$  \;
  \caption{Hybrid  transport method \label{alg:hybrid_algorithm}}
\end{algorithm}

\section{Numerical Results and Analysis}

The errors of  numerical solutions of the HQD and HSM methods are analyzed on the transport problem in the spatial domain $0 \le x  \le 1$ with  $\Sigma_t = 1.0$,  $\Sigma_s = 0.9$,   $Q = 1.0$, and vacuum BCs.
The spatial mesh  with $I$ cells is uniform.
The accuracy of numerical solutions is evaluated by means of the reference deterministic solution computed on a sequence of refined grids in the phase space with Aitken extrapolation. 
The obtained benchmark  solution  ($\phi^{ex}$) is exact in its 6 significant digits.
We present the relative error in $L_2$-norm given by

\begin{equation} \label{eqn:rel_L2_Norm_Error}
  \text{RE}_{L_2} = \sqrt{ \frac{\sum_{i=1}^I (\phi_i - \phi_i^{ex})^2\Delta x_i}
  {\sum_{i=1}^I (\phi_i^{ex})^2 \Delta x_i} } \, .
\end{equation}

We performed 100 calculations of the test by the HQD and HSM methods and compare the error of each hybrid solution with one of the analog MC solution the particle histories of which were used to compute functionals for that hybrid algorithm calculation.
Tables \ref{tab:loqd_to_mc} and \ref{tab:losm_to_mc} present the fraction of times that each hybrid method produced the solution with a lower relative $L_2$ error compared to the companion MC calculations.
These results show the combination of grids and number of histories where discretization error dominates in the error in the solution of the HQD and HSM methods.
The hybrid methods produced on average a more accurate solution than direct MC calculations for most of considered combination of $N$ and $\Delta x$'s.
HSM has a higher rate of wins vs. MC compared to the HQD method.
We also compared $L_\infty$ norms of error as well, though not presented in this summary, which shows similar behavior to the integral norm results.
\begin{table}[th]
  \centering
  \caption{HQD $L_2$ win ratio compared to MC as a function of $N$ particles and cell width $\Delta x$ }
  {
  \begin{tabular}{|c|c|c|c|c|}
      \hline
      $I$   & $\Delta x$ & $N=10^3$ & $N=10^4$ & $N=10^5$ \\
      \hline
      4& $2^{-2}$ & 0.55 & 0.13 &0.00 \\
      8& $2^{-3}$ & 0.71 & 0.63 &0.17 \\
     16& $2^{-4}$ & 0.74 & 0.72 &0.89 \\
     32& $2^{-5}$ & 0.70 & 0.68 &0.70 \\
     64& $2^{-6}$ & 0.73 & 0.77 &0.77 \\
     \hline
  \end{tabular}
  }
  \label{tab:loqd_to_mc}
\end{table}
\begin{table}[th]
  \centering
  \caption{HSM $L_2$ win ratio compared to MC as a function of $N$ particles and cell width $\Delta x$}
  {\small
  \begin{tabular}{|c|c|c|c|c|}
      \hline
     $I$ & $\Delta x$  & $N=10^3$ & $N=10^4$ & $N=10^5$ \\
      \hline
      4& $2^{-2}$ & 0.66 & 0.22 & 0.00\\
      8& $2^{-3}$ & 0.78 & 0.67 & 0.29\\
     16& $2^{-4}$ & 0.80 & 0.76 & 0.79\\
     32& $2^{-5}$ & 0.82 & 0.73 & 0.80\\
     64& $2^{-6}$ & 0.84 & 0.79 & 0.76\\
     \hline
  \end{tabular}
  }
  \label{tab:losm_to_mc}
\end{table}

\begin{figure*}[tb]
  \centering
  \begin{subfigure}[b]{0.49\textwidth}
      \centering
      \includegraphics[width=0.9\textwidth]{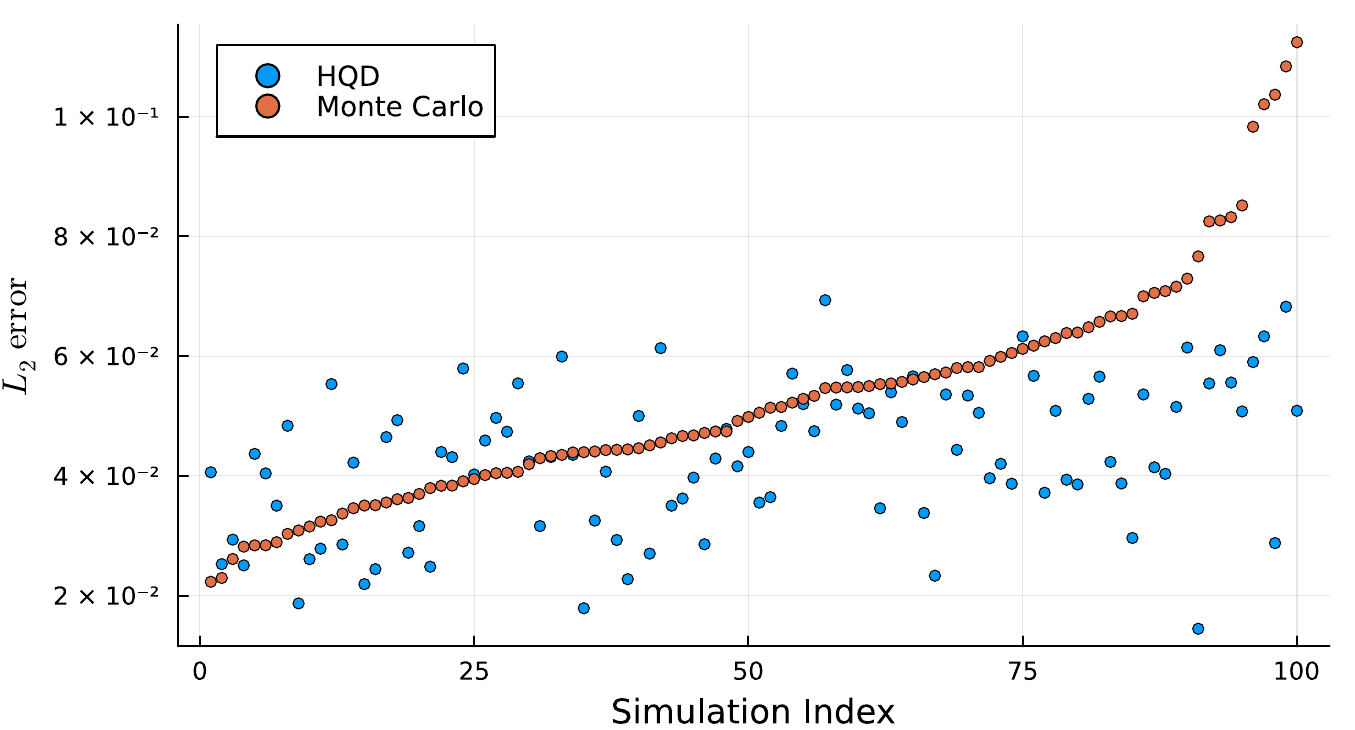}
        \caption*{$N = 10^3$ histories, relative errors  in $L_2$-norm of MC and HQD}
  \end{subfigure}
  \begin{subfigure}[b]{0.49\textwidth}
      \centering
      \includegraphics[width=0.9\textwidth]{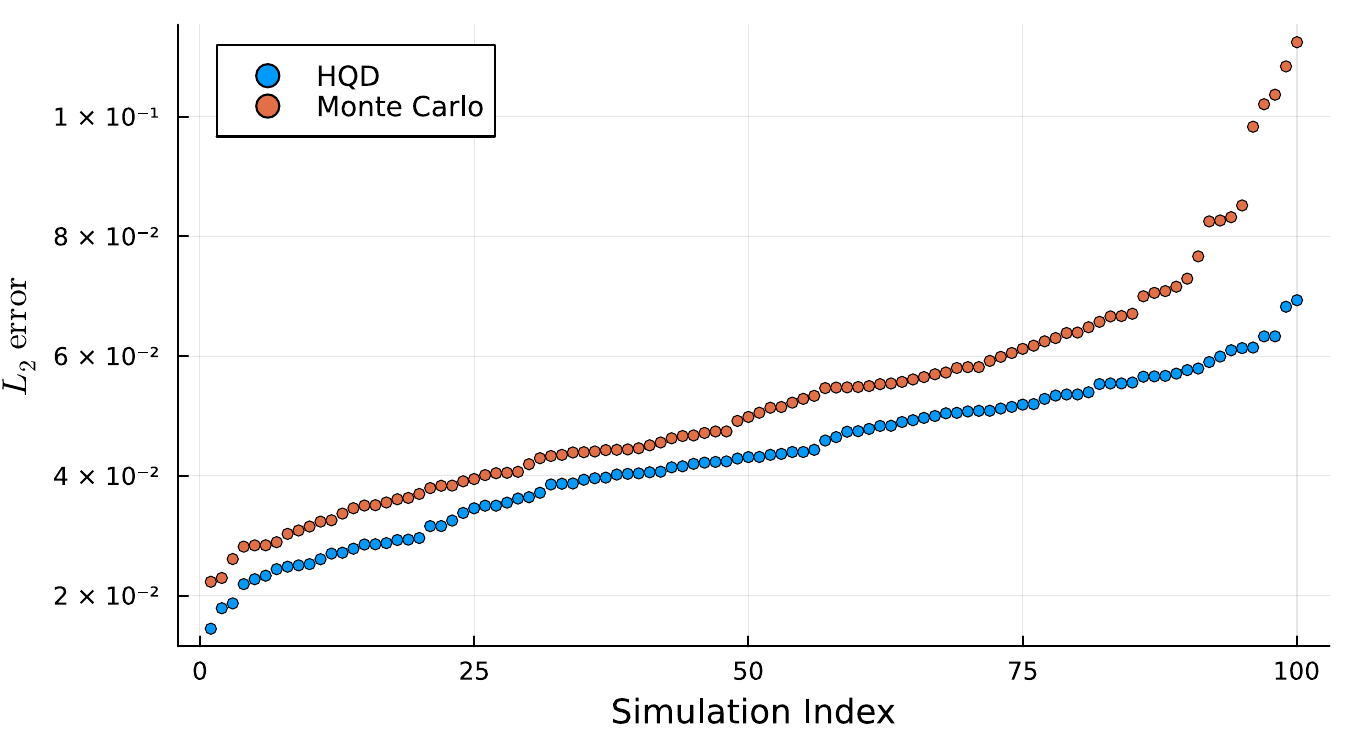}
     \caption*{$N = 10^3$ histories, sorted relative error  in $L_2$-norm of MC and HQD}
  \end{subfigure}
  \begin{subfigure}[b]{0.49\textwidth}
      \centering
      \includegraphics[width=0.9\textwidth]{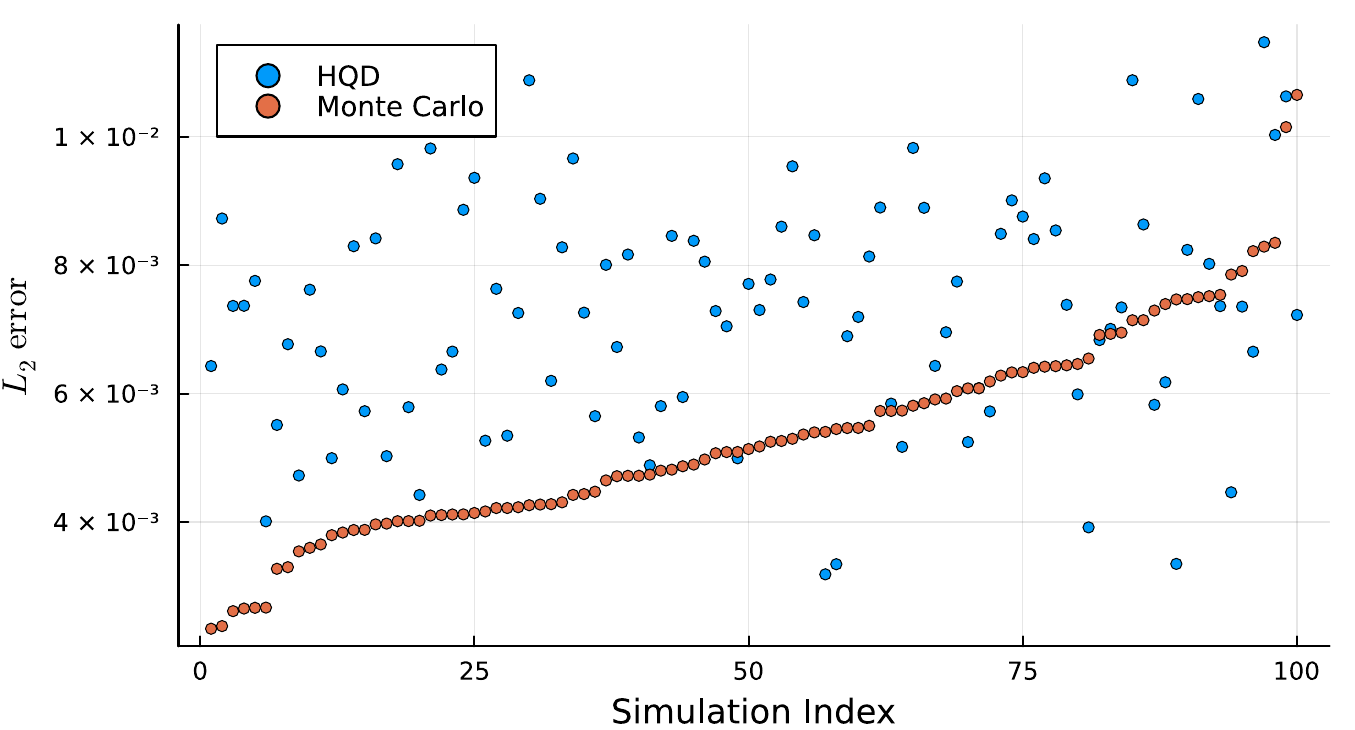}
       \caption*{$N = 10^5$ histories, relative errors  in $L_2$-norm of MC and HQD}
  \end{subfigure}
  \begin{subfigure}[b]{0.49\textwidth}
      \centering
      \includegraphics[width=0.9\textwidth]{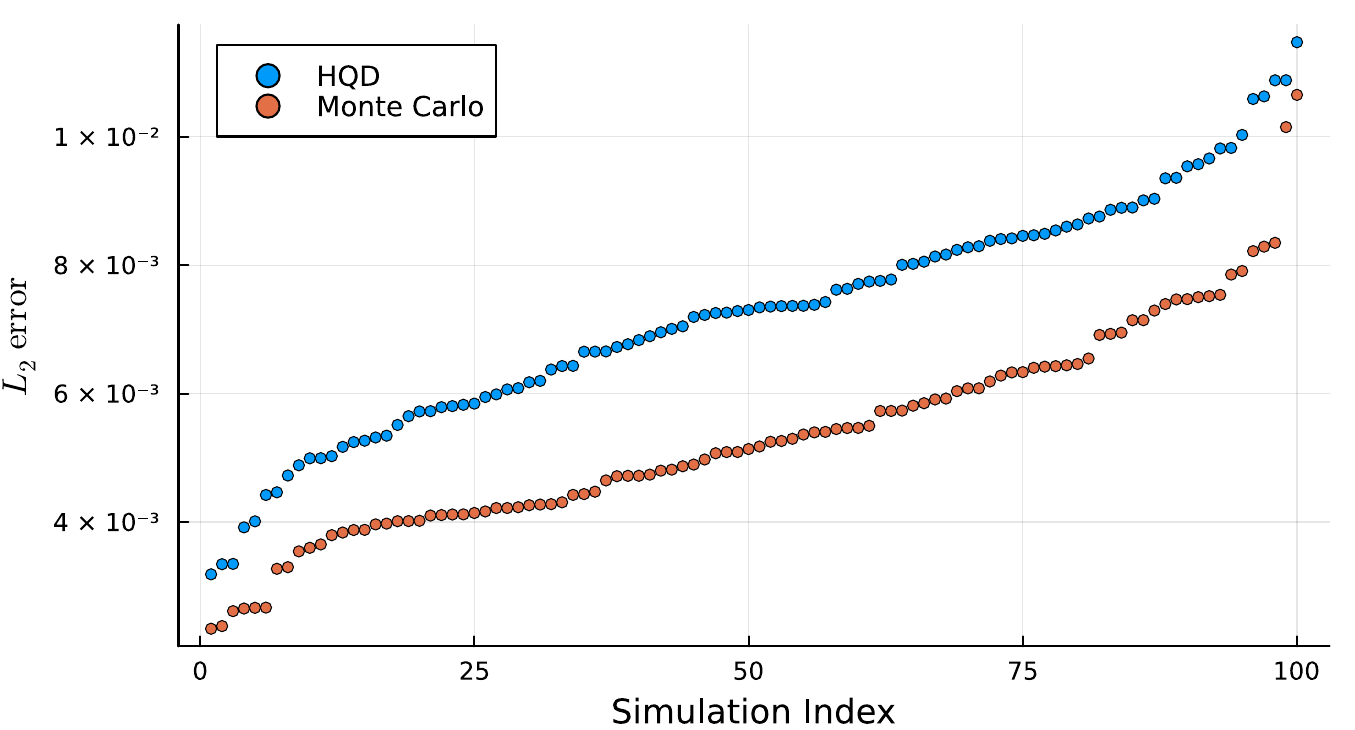}
        \caption*{$N = 10^5$ histories, sorted  relative error  in $L_2$-norm of MC and HQD}
  \end{subfigure}
  \caption{Relative error in  $L_2$-norm of MC and HQD solutions, I = 8, $\Delta x= 2^{-3}$.}
  \label{fig:8_cell_line}
\end{figure*}
\begin{figure*}[tb]
  \centering
  \begin{subfigure}[b]{0.49\textwidth}
      \centering
      \includegraphics[width=0.9\textwidth]{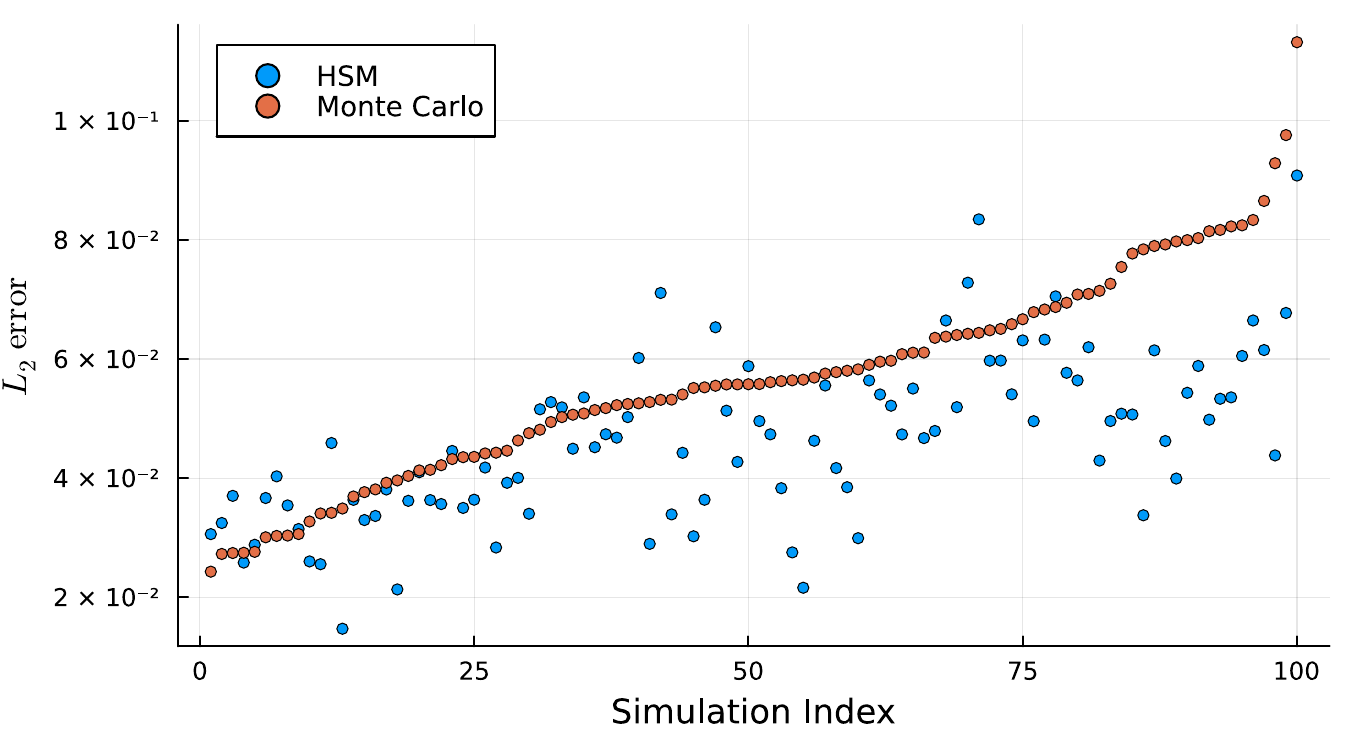}
        \caption*{$N = 10^3$ histories, relative errors  in $L_2$-norm of MC and HSM}
  \end{subfigure}
  \begin{subfigure}[b]{0.49\textwidth}
      \centering
      \includegraphics[width=0.9\textwidth]{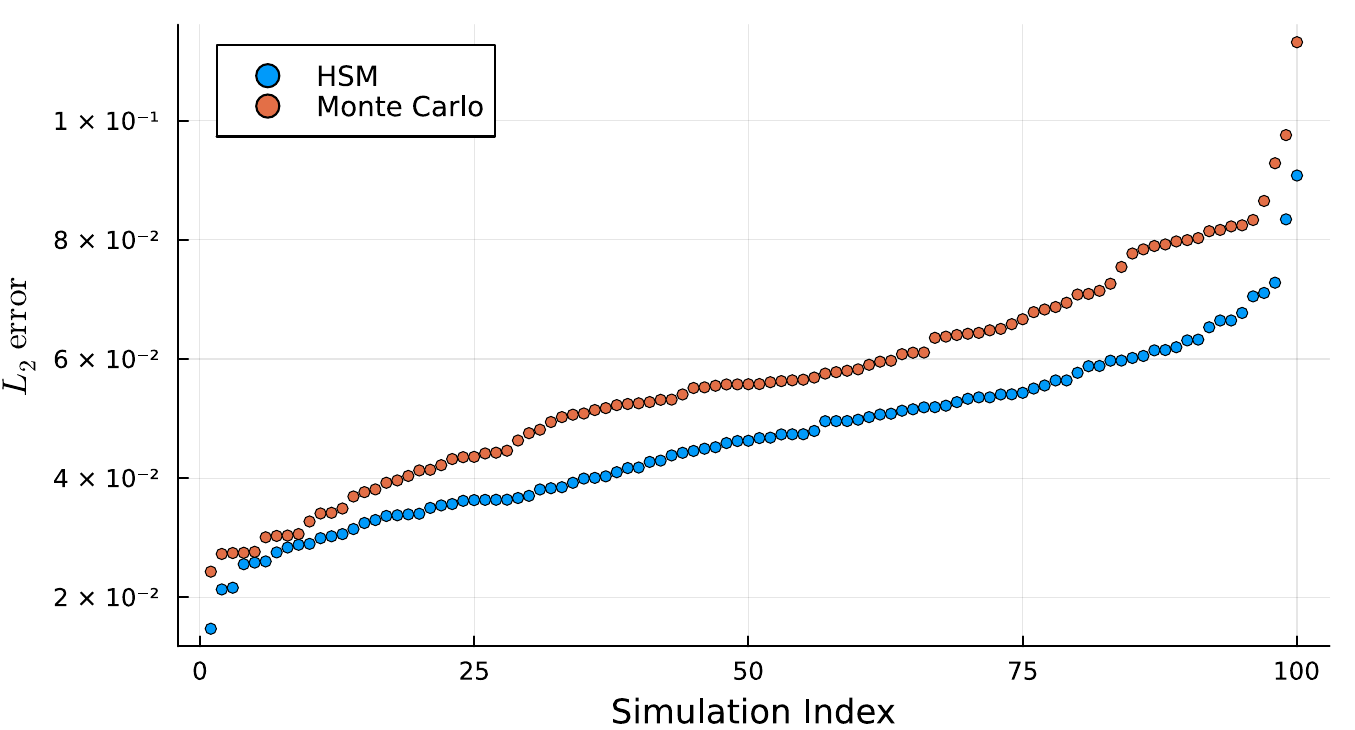}
     \caption*{$N = 10^3$ histories, sorted relative error  in $L_2$-norm of MC and HSM}
  \end{subfigure}
  \begin{subfigure}[b]{0.49\textwidth}
      \centering
      \includegraphics[width=0.9\textwidth]{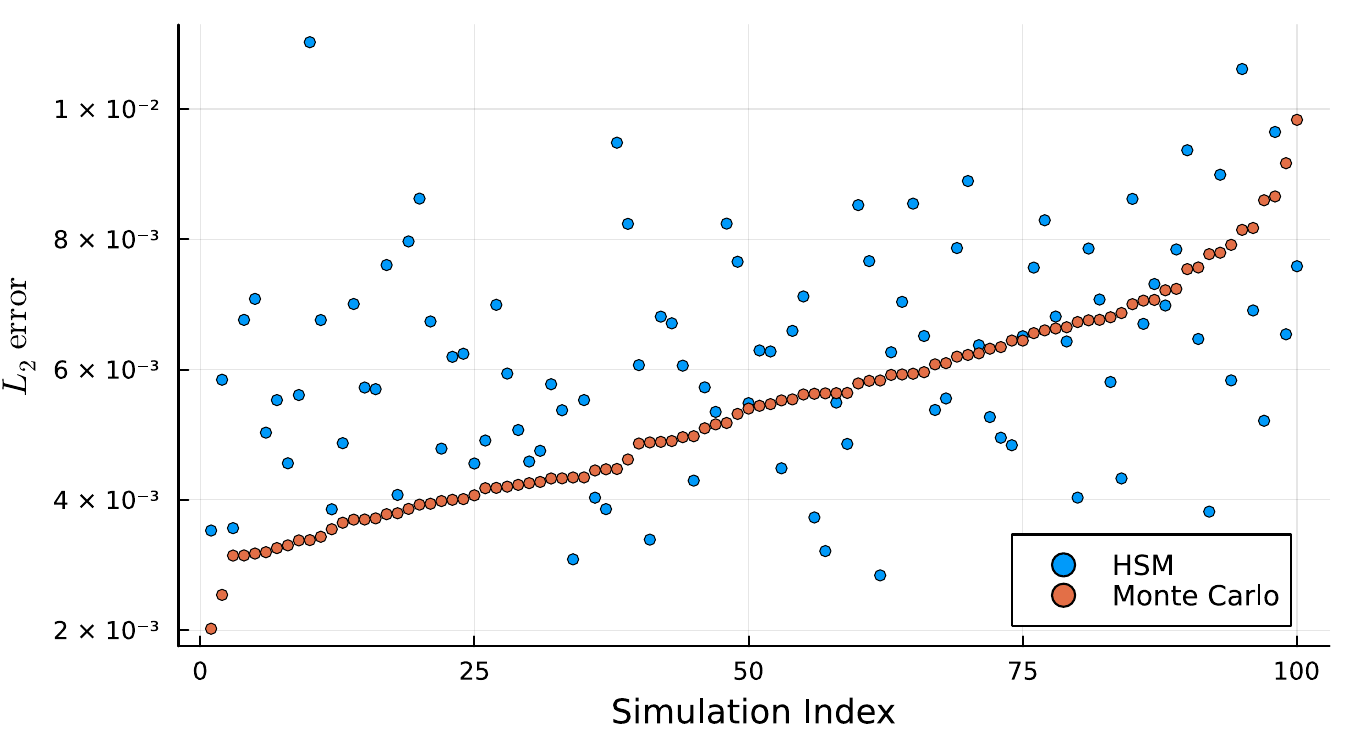}
       \caption*{$N = 10^5$ histories, relative errors  in $L_2$-norm of MC and HSM}
  \end{subfigure}
  \begin{subfigure}[b]{0.49\textwidth}
      \centering
      \includegraphics[width=0.9\textwidth]{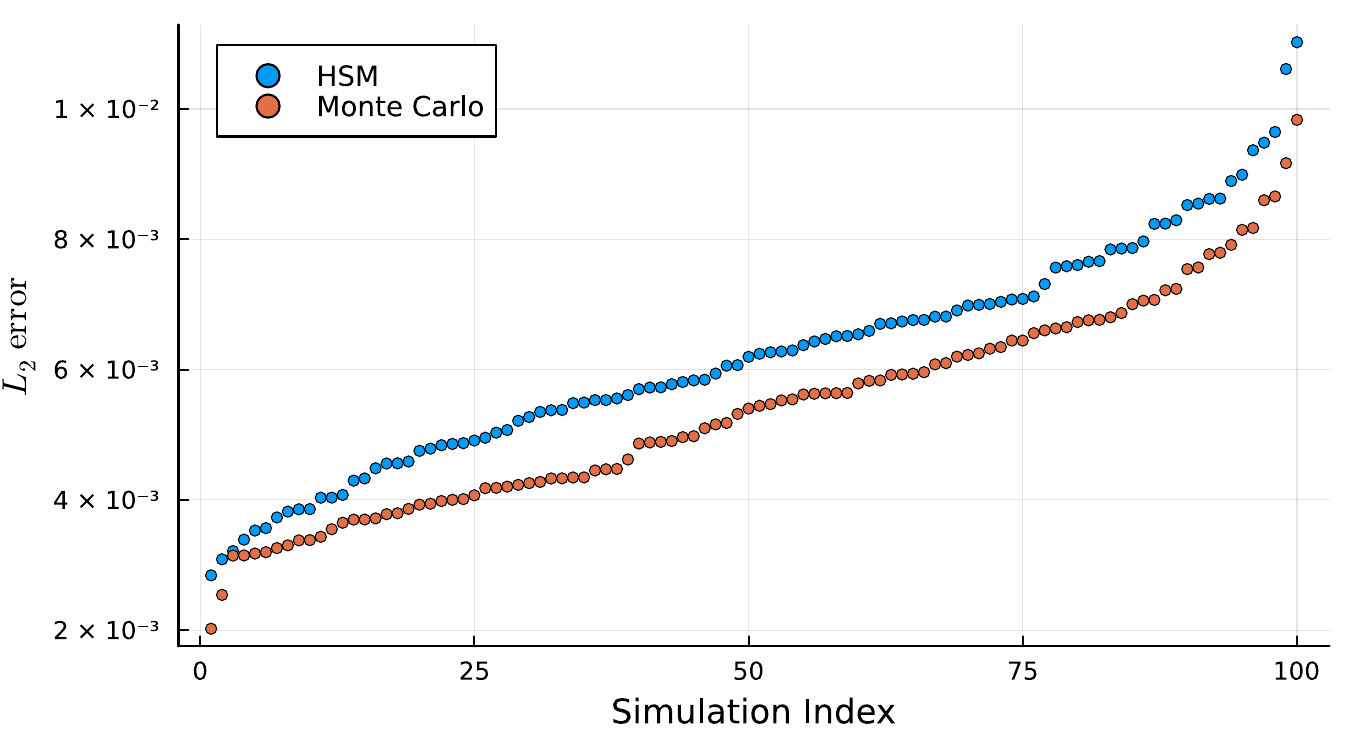}
        \caption*{$N = 10^5$ histories, sorted  relative error  in $L_2$-norm of MC and HSM}
  \end{subfigure}
  \caption{Relative error in  $L_2$-norm of MC and HSM solutions, I = 8, $\Delta x= 2^{-3}$.}
  \label{fig:8_cell_line_sm}
\end{figure*}

\begin{figure*}[tb]
  \centering
  \begin{subfigure}[b]{0.49\textwidth}
      \centering
      \includegraphics[width=0.9\textwidth]{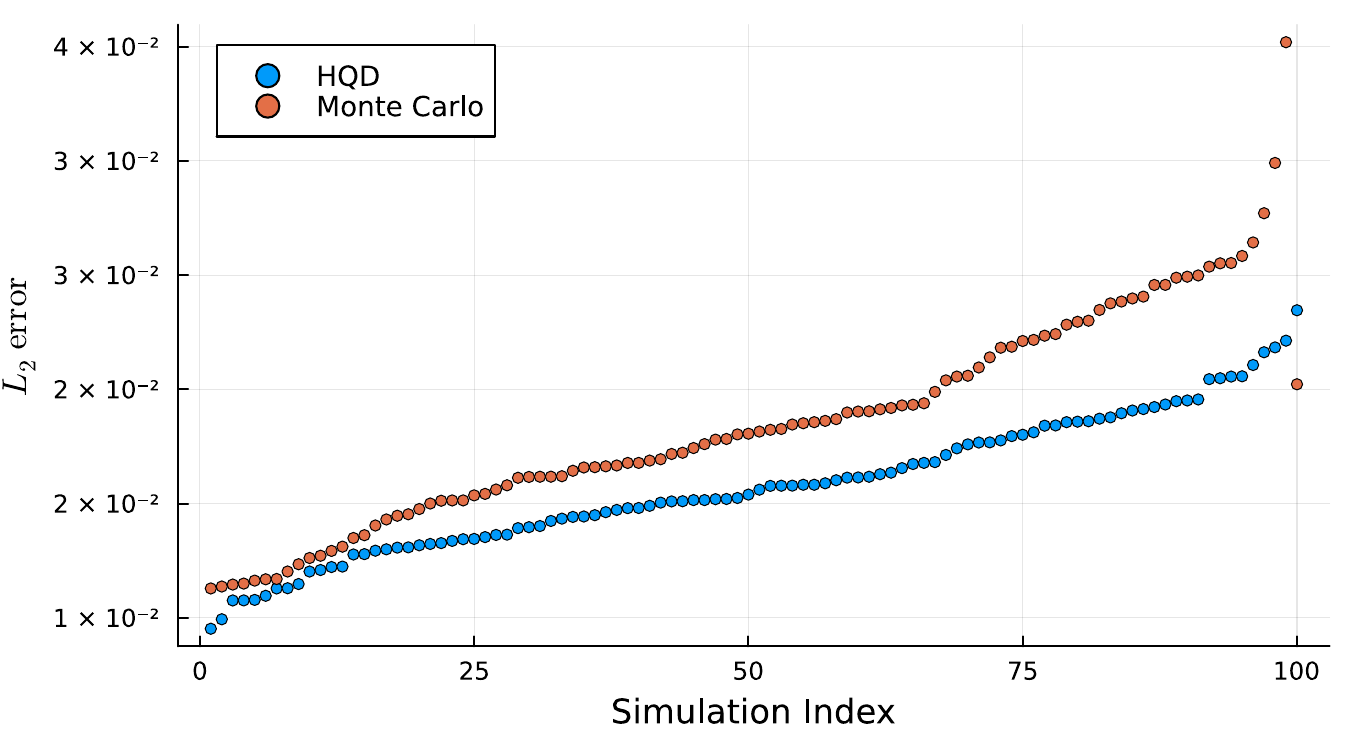}
        \caption*{$N = 10^4$ histories, MC and HQD}
  \end{subfigure}
  \begin{subfigure}[b]{0.49\textwidth}
      \centering
      \includegraphics[width=0.9\textwidth]{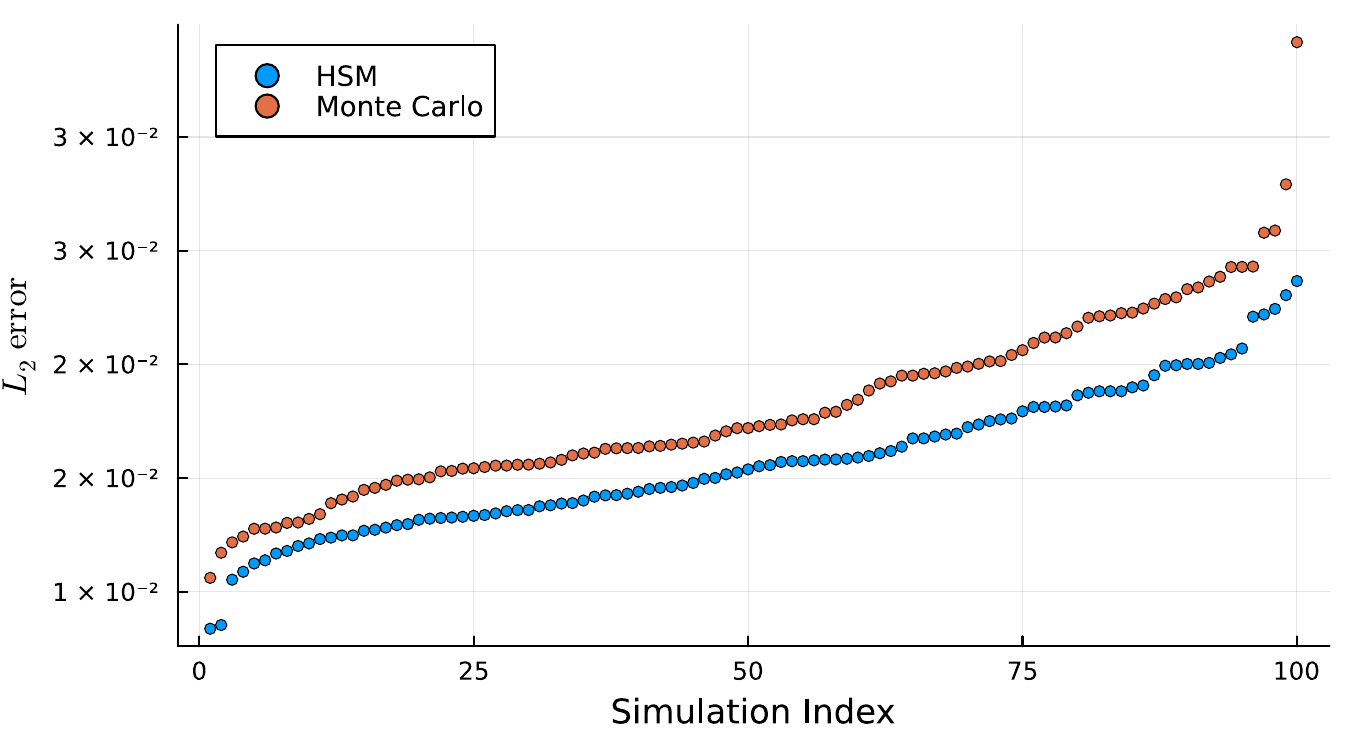}
        \caption*{$N = 10^4$ histories, MC and HSM}
  \end{subfigure}
  \begin{subfigure}[b]{0.49\textwidth}
      \centering
      \includegraphics[width=0.9\textwidth]{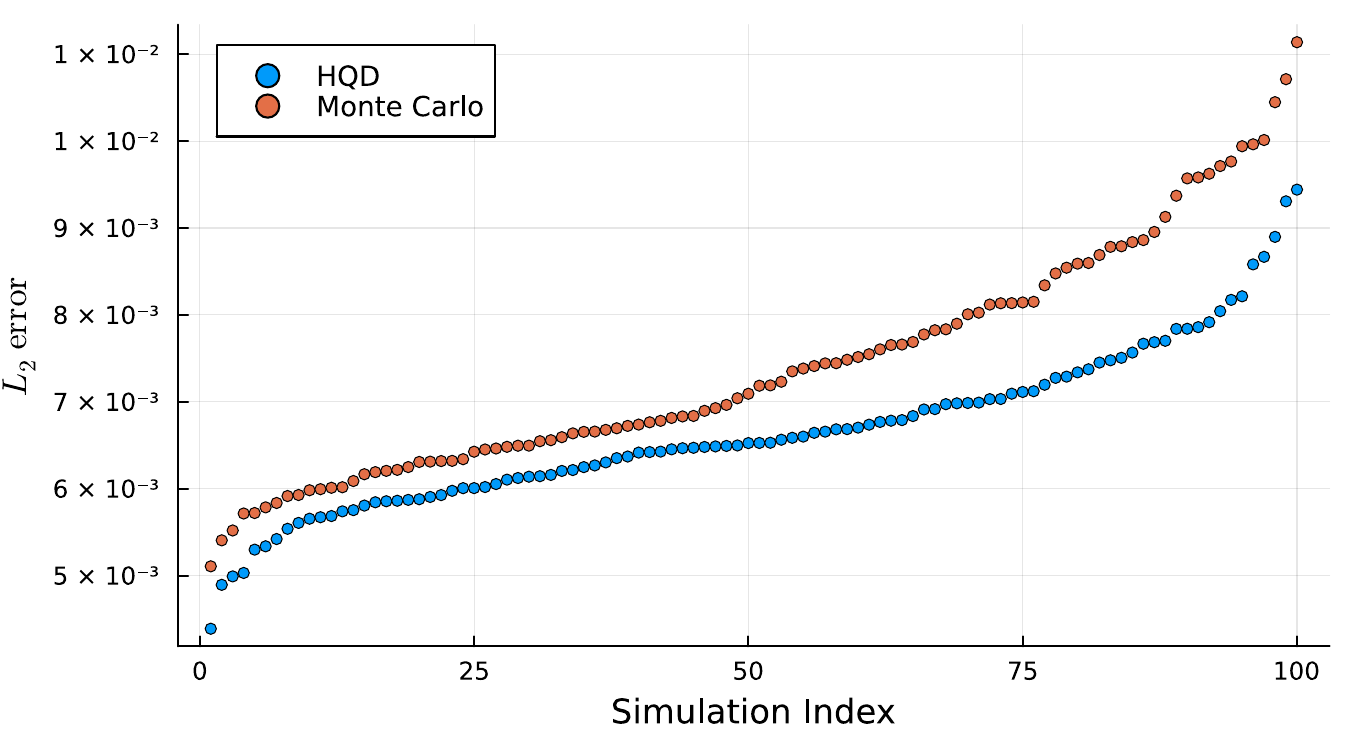}
        \caption*{$N = 10^5$ histories, MC and HQD}
  \end{subfigure}
  \begin{subfigure}[b]{0.49\textwidth}
      \centering
      \includegraphics[width=0.9\textwidth]{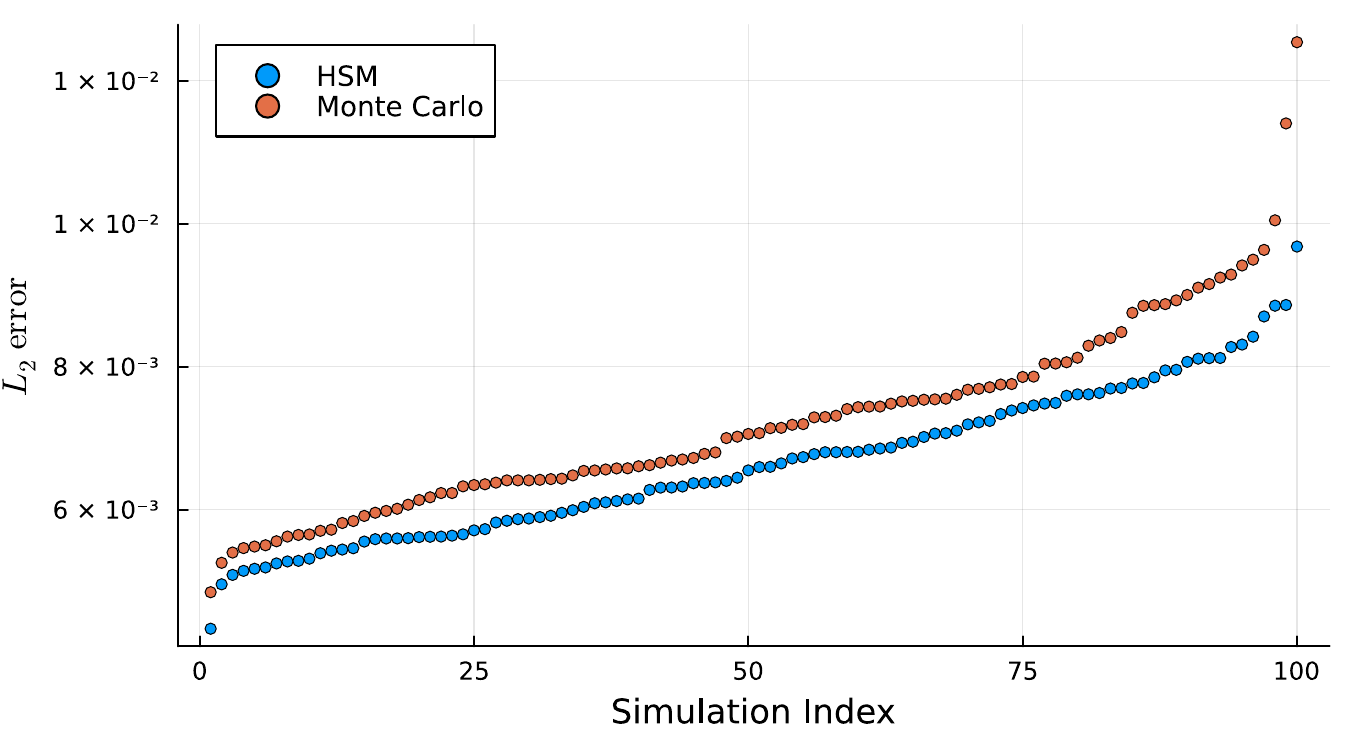}
        \caption*{$N = 10^5$ histories, MC and HSM}
  \end{subfigure}
    \caption{Sorted    relative error in  $L_2$-norm of MC and hybrid solutions, I = 64, $\Delta x= 2^{-6}$.}
  \label{fig:64_cell_line}
\end{figure*}

\pagebreak
Figure \ref{fig:8_cell_line} shows the  relative $L_2$ errors of MC and HQD for different number of particle histories, namely, $N=10^3$ and $N=10^5$.
On the left side of Figure \ref{fig:8_cell_line} we present  the  errors of MC and the corresponding HQD solutions computed with the same ensemble of particle histories versus the index of simulation. 
These plots illustrate the performance of  two  methods compared to each other in every case.
The plots of sorted $L_2$ errors of both MC and HQD are shown on the right sub-figure.
These results illustrate the effect of increasing the number of histories on a coarse spatial grid.
An observation from the unsorted results is the generally weak correlation between the MC and HQD errors.
Our sorted $L_2$ results demonstrate the average behavior of the MC and HQD methods showing which method tends to have lower relative $L_2$ error.
We notice that for low number of particle histories the HQD methods performs better than MC, especially for cases with high MC statistical error.
The discretization error becomes dominant  with  increasing  the particle number used in calculations leading to less accurate solutions compared to MC (see Fig. \ref{fig:8_cell_line} for $N=10^5$).
Figure \ref{fig:8_cell_line_sm} presents the same set plots for the HSM method with MC simulations. The performance  of the HSM  method is similar one of the HQD method.
Figure  \ref{fig:64_cell_line}   shows  the sorted errors for   fine mesh and large number particle history cases for both HQD and HSM methods.
We see as a trend that the hybrid methods reduce error in the simulation and eliminate the right tail, i.e. sharp increase in error, present in the Monte Carlo results.
We caution the reader to interpret Figure  \ref{fig:64_cell_line}  in light of the win ratio shown in Tables \ref{tab:loqd_to_mc} and \ref{tab:losm_to_mc}.

\begin{figure*}[htb]
\centering
  \begin{subfigure}[b]{0.5\textwidth}
      \centering
      \includegraphics[width=\textwidth]{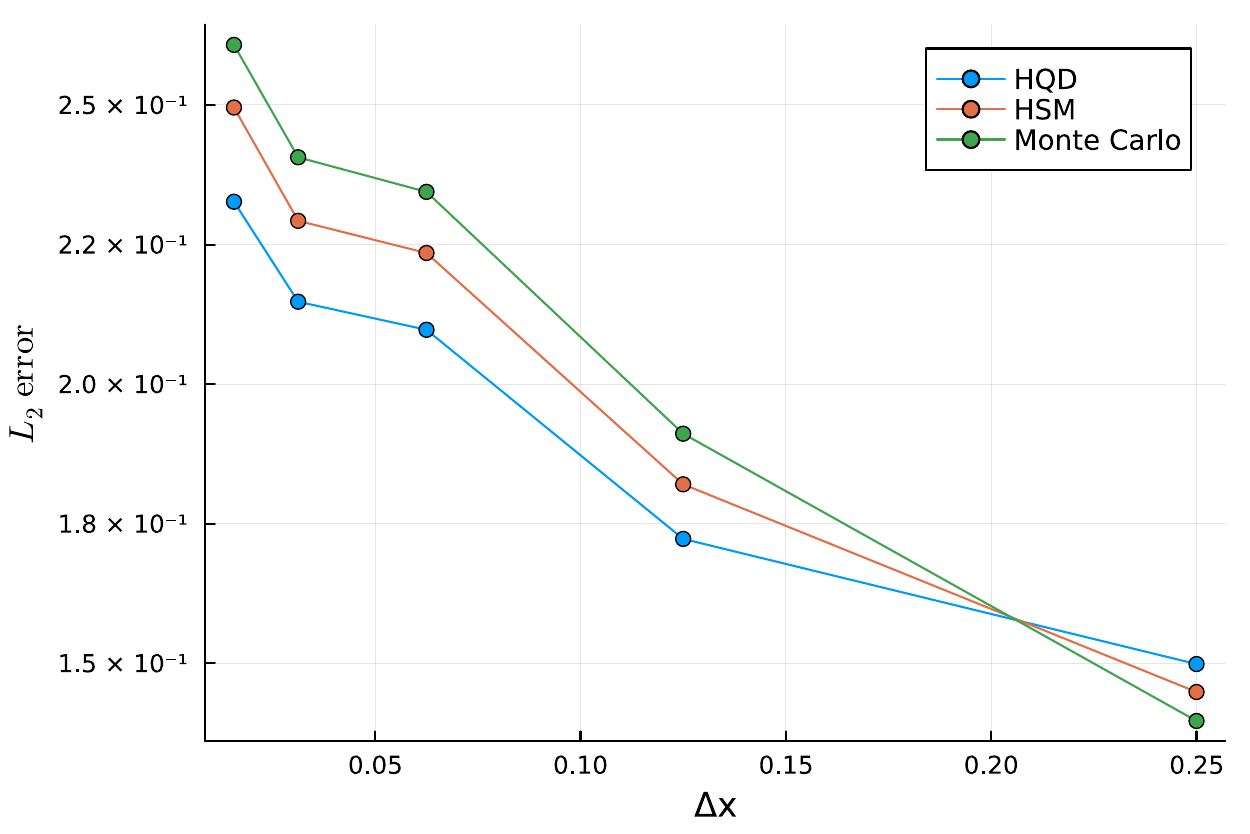}
      \caption*{$N = 10^2$, same histories for each method}
  \end{subfigure}
  \begin{subfigure}[b]{0.5\textwidth}
      \centering
      \includegraphics[width=\textwidth]{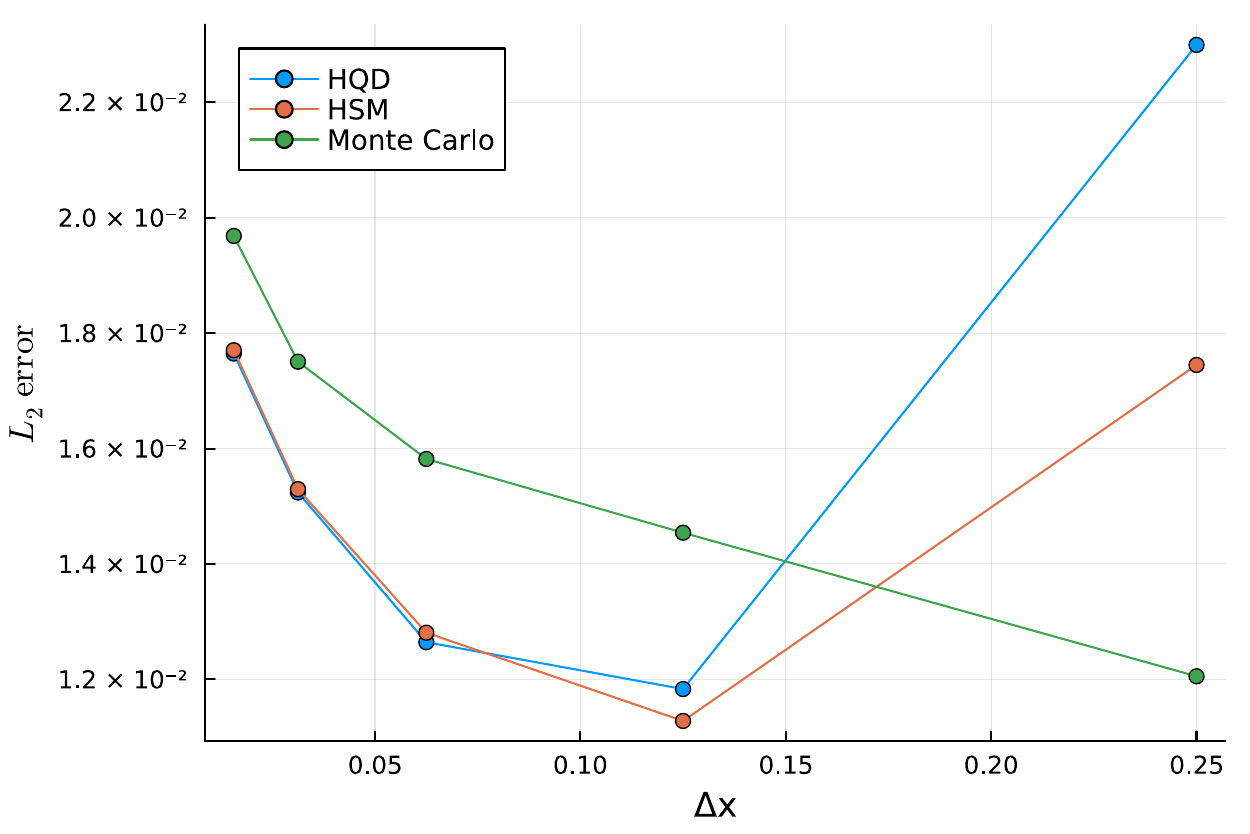}
      \caption*{$N = 10^4$, same histories for each method}
  \end{subfigure}
  \begin{subfigure}[b]{0.5\textwidth}
      \centering
      \includegraphics[width=\textwidth]{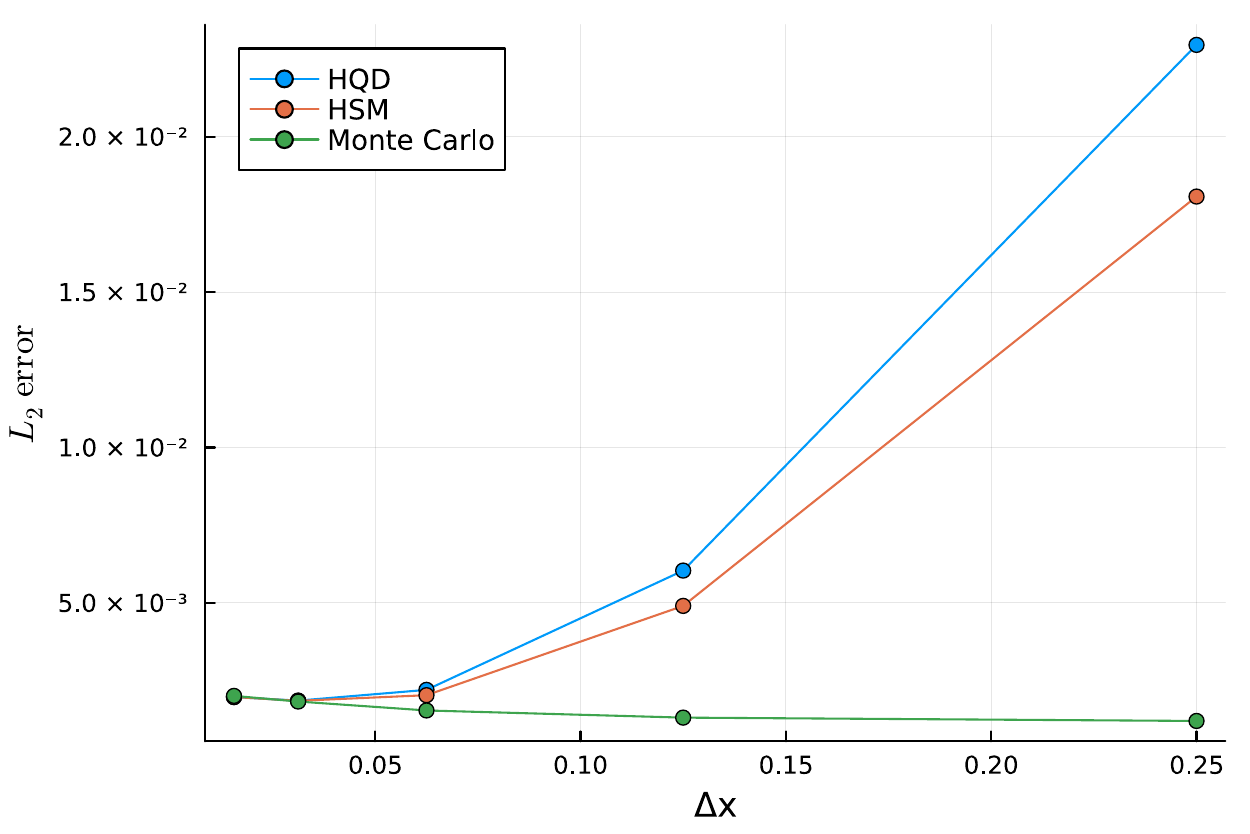}
      \caption*{$N = 10^6$, same histories for each method}
  \end{subfigure}
  \caption{Relative error in  $L_2$-norm of MC and hybrid solutions vs. $\Delta x$
  under grid refinement for varying number of  histories}
  \label{fig:hybrid_mc_l2_grid_ref}
\end{figure*}

Table \ref{tab:mc_hybrid_l2_errors} presents the relative error in $L_2$-norm ($RE_{L_2}$) of the MC solution with implicit capture, and hybrid solutions for the set of calculations with spatial grid refinement and increasing number of particle histories.
In the cases $N=10^2$ and $N=10^4$, the HQD and HSM solutions are slightly more accurate on all grids except for $\Delta x = 2^{-2}$.
Note that the error of MC solution increases under grid refinement, since as the tally cells shrinks the total tracks-length in a cell decreases leading to an increase in uncertainty.
The HQD and HSM methods show an increase in error under grid refinement for low particle counts ($N=10^2$) due to the large statistical uncertainty.
For large number of histories ($N=10^6$), we notice significant reduction in error of hybrid solutions as $\Delta x$ decreases driven by the spatial convergence of the underlying FV schemes for LOQD and LOSM equations.
Table \ref{tab:mc_hybrid_l2_errors_ratio} presents the ratios of errors in $L_2$-norms $RE_{L_2}(2 \Delta x)/RE_{L_2}(\Delta x)$ for the solution on neighboring grids.
These results quantify the change of error in solution with grid refinement.
Figure \ref{fig:hybrid_mc_l2_grid_ref} shows the relative error in $L_2$-norm vs. $\Delta x$ for each method for a select number of particle histories.
The results demonstrate how hybrid solutions are more accurate for small ($N=10^2$) and sufficiently large number ($N=10^4$) of histories.
We also notice that for very high particle counts ($N=10^6$) discretization error is dominant for these hybrid methods with FV spatial schemes.

\begin{table}[htb]
  \centering
  \caption{Relative error in  $L_2$-norm ($RE_{L_2}$) for MC, HQD, and HSM on sequence of grids}
    \begin{tabular}{|c|c|c|c|}
        \hline
        $N = 10^2$         &  MC                   & HQD                   & HSM                   \\  
        \hline
        $\Delta x =2^{-2}$ & 1.4 $\times  10^{-1}$ & 1.5$ \times  10^{-1}$ & 1.45$\times  10^{-1}$ \\  
        $\Delta x =2^{-3}$ & 1.91$\times  10^{-1}$ & 1.72$\times  10^{-1}$ & 1.82$\times  10^{-1}$ \\ 
        $\Delta x =2^{-4}$ & 2.34$\times  10^{-1}$ & 2.10$\times  10^{-1}$ & 2.23$\times  10^{-1}$ \\  
        $\Delta x =2^{-5}$ & 2.41$\times  10^{-1}$ & 2.15$\times  10^{-1}$ & 2.29$\times  10^{-1}$ \\ 
        $\Delta x =2^{-6}$ & 2.61$\times  10^{-1}$ & 2.33$\times  10^{-1}$ & 2.50$\times  10^{-1}$ \\ 
        \hline  \hline
        $N = 10^4$         &  MC                   & HQD                   & HSM                   \\  
        \hline
        $\Delta x =2^{-2}$ & 1.21$\times  10^{-2}$ & 2.30$\times  10^{-2}$ & 1.74$\times  10^{-2}$ \\  
        $\Delta x =2^{-3}$ & 1.45$\times  10^{-2}$ & 1.18$\times  10^{-2}$ & 1.13$\times  10^{-2}$ \\  
        $\Delta x =2^{-4}$ & 1.58$\times  10^{-2}$ & 1.26$\times  10^{-2}$ & 1.28$\times  10^{-2}$ \\  
        $\Delta x =2^{-5}$ & 1.75$\times  10^{-2}$ & 1.52$\times  10^{-2}$ & 1.53$\times  10^{-2}$ \\  
        $\Delta x =2^{-6}$ & 1.97$\times  10^{-2}$ & 1.76$\times  10^{-2}$ & 1.77$\times  10^{-2}$ \\  
        \hline  \hline
        $N = 10^6$         &  MC                   & HQD                   & HSM                   \\  
        \hline
        $\Delta x =2^{-2}$ & 1.20$\times  10^{-3}$ & 2.30$\times  10^{-2}$ & 1.81$\times  10^{-2}$ \\  
        $\Delta x =2^{-3}$ & 1.31$\times  10^{-3}$ & 6.04$\times  10^{-3}$ & 4.90$\times  10^{-3}$ \\  
        $\Delta x =2^{-4}$ & 1.54$\times  10^{-3}$ & 2.20$\times  10^{-3}$ & 2.03$\times  10^{-3}$ \\  
        $\Delta x =2^{-5}$ & 1.83$\times  10^{-3}$ & 1.86$\times  10^{-3}$ & 1.85$\times  10^{-3}$ \\  
        $\Delta x =2^{-6}$ & 2.01$\times  10^{-3}$ & 1.97$\times  10^{-3}$ & 1.97$\times  10^{-3}$ \\  
        \hline
    \end{tabular}
  \label{tab:mc_hybrid_l2_errors}
\end{table}

\begin{table}[htb]
  \centering
  \caption{Ratios of errors in $L_2$-norms, $RE_{L_2}(2 \Delta x)/RE_{L_2}(\Delta x)$}
    \begin{tabular}{|c|c|c|c|}
        \hline
        $N = 10^2$         & MC   & HQD  & HSM  \\ 
        \hline
        \hline
        $\Delta x =2^{-3}$ & 0.73 & 0.87 & 0.80 \\ 
        $\Delta x =2^{-4}$ & 0.82 & 0.82 & 0.82 \\  
        $\Delta x =2^{-5}$ & 0.97 & 0.98 & 0.98 \\  
        $\Delta x =2^{-6}$ & 0.92 & 0.92 & 0.92 \\  
        \hline  
        \hline
        $N = 10^4$         & MC   & HQD  & HSM  \\  
        \hline
        $\Delta x =2^{-3}$ & 0.83 & 1.94 & 1.55 \\  
        $\Delta x =2^{-4}$ & 0.92 & 0.94 & 0.88 \\  
        $\Delta x =2^{-5}$ & 0.90 & 0.83 & 0.84 \\ 
        $\Delta x =2^{-6}$ & 0.89 & 0.86 & 0.86 \\  
        \hline 
        \hline
        $N = 10^6$         & MC   & HQD  & HSM  \\  
        \hline
        $\Delta x =2^{-3}$ & 0.92 & 3.80 & 3.69 \\  
        $\Delta x =2^{-4}$ & 0.85 & 2.74 & 2.42 \\ 
        $\Delta x =2^{-5}$ & 0.84 & 1.19 & 1.10 \\ 
        $\Delta x =2^{-6}$ & 0.91 & 0.94 & 0.94 \\  
        \hline
    \end{tabular}
    \label{tab:mc_hybrid_l2_errors_ratio}
\end{table}

\section{Conclusions}
The HQD and HSM methods with the FV approximation of the low-order equations with second-order convergence generate more accurate numerical solutions than MC on average when discretization error is small.
These hybrid techniques also showed advantages in simulations with low particle counts per tally cell.
For multi-dimensional problems, it is expensive to obtain large number of particle histories for each phase-space elements.
Therefore, average particle history per phase space element should be considered when evaluating low-order hybrid solvers.
The HSM method shows some advantages when compared to HQD.
Further testing and analysis are needed to study effects of statistical noise and discretization error. 

 \clearpage
\section{Acknowledgments}
\vspace{-0.1cm}
The work of the first author (VNN) was supported under an University Nuclear Leadership Program Graduate Fellowship. Any opinions, findings, conclusions or recommendations expressed in this publication are those of the author(s) and do not necessarily reflect the views of the Department of Energy Office of Nuclear Energy. 
The work of the second author (DYA) was supported by the Center for Exascale Monte-Carlo Neutron Transport (CEMeNT), a PSAAP-III project funded by the Department of Energy, grant number DE-NA003967.

\bibliographystyle{elsarticle-num}
\bibliography{bibliography-VN-DYA-2024}

\begin{thebibliography}{1}
\expandafter\ifx\csname url\endcsname\relax
  \def\url#1{\texttt{#1}}\fi
\expandafter\ifx\csname urlprefix\endcsname\relax\def\urlprefix{URL }\fi
\expandafter\ifx\csname href\endcsname\relax
  \def\href#1#2{#2} \def\path#1{#1}\fi

\bibitem{larsen-yang-nse-2008}
{E. W. Larsen}, {J. Yang}, A functional {M}onte {C}arlo method for k-eigenvalue
  problems, Nuclear Science and Engineering 159 (2008) 107--126.

\bibitem{lee-physor-2010}
{M.J. Lee}, {H.G. Joo}, {D. Lee}, {K. Smith}, Investigation of {CMFD}
  accelerated {M}onte {C}arlo eigenvalue calculation with simplified low
  dimensional multigroup formulation, in: Proc. of PHYSOR 2010, Int. Conf. on
  the Physics of Reactors, Pittsburgh, PA, May 9-14, 2010.

\bibitem{wolters-nse-2013}
{E. R. Wolters}, {E. W. Larsen}, {W. R. Martin}, Hybrid {M}onte
  {C}arlo–{CMFD} methods for accelerating fission source convergence, Nuclear
  Science and Engineering 174 (2013) 286–299.

\bibitem{pozulp-mc2023}
M.~Pozulp, T.~Haut, P.~Brantley, J.~Vujic, An implicit monte carlo acceleration
  scheme, in: Proc. of M\&C 2023, Int. Conf. on Math. \& Comp. Methods Applied
  to Nucl. Sci \& Eng., Niagara Falls, Canada, August 13-17, 2023.

\bibitem{gol'din-cmmp-1964}
{V. Ya. Gol'din}, A quasi-diffusion method of solving the kinetic equation,
  Comp. Math. and Math. Phys. 4 (1964) 136--149.

\bibitem{sm-1976}
{ E. E. Lewis}, {W. F. Miller, Jr.}, A comparison of p1 synthetic acceleration
  techniques, Transactions of the American Nuclear Society 23 (1976) 202--203.

\bibitem{mla-ewl-pne-2002}
{M. L. Adams}, {E. W. Larsen}, Fast iterative methods for discrete-ordinates
  particle transport calculations, Progress in Nuclear Energy 40 (2002) 3--159.

\end{thebibliography}

\end{document}